\documentclass[11pt]{article}
\usepackage{latexsym,amssymb}
\usepackage{amsmath}
\usepackage[mathscr]{eucal}
\usepackage{color}

\newtheorem{Theorem}{\bf Theorem}[section]
\newtheorem{Lemma}{\bf Lemma}[section]
\newtheorem{Proposition}{\bf Proposition}[section]  
\newtheorem{Corollary}{\bf Corollary}[section]
\newtheorem{Remark}{\bf Remark}[section]

\newtheorem{Definition}{\bf Definition}[section]

\newenvironment{theorem}{\begin{Theorem}$\!\!\!$}{\end{Theorem}}
\newenvironment{lemma}{\begin{Lemma}$\!\!\!$}{\end{Lemma}}

\newenvironment{corollary}{\begin{Corollary}$\!\!\!$}{\end{Corollary}}
\newenvironment{remark}{\begin{Remark}$\!\!\!$}{\end{Remark}}
\newenvironment{definition}{\begin{Definition}$\!\!\!$}{\end{Definition}}

\numberwithin{equation}{section}

\begin{document}
\title{To logconcavity and beyond}
\author{Kazuhiro Ishige, Paolo Salani and Asuka Takatsu}

\date{}
\maketitle
\begin{abstract}
In 1976 Brascamp and Lieb proved that 
the heat flow preserves logconcavity. 
In this paper, introducing a variation of concavity, 
we show that it preserves in fact a stronger property than logconcavity 
and we identify the strongest concavity preserved by the heat flow. 

\end{abstract}
\vspace{40pt}
\noindent Addresses:

\smallskip
\noindent 
K. I.: Graduate School of Mathematical Sciences, The University of Tokyo, 3-8-1 
Komaba, Meguro-ku, Tokyo 153-8914, Japan.\\
\noindent 
E-mail: {\tt ishige@ms.u-tokyo.ac.jp}\\

\smallskip
\noindent 
P. S.: Dipartimento di Matematica ``U. Dini'', 
Universit\`a di Firenze, viale Morgagni 67/A, 50134 Firenze\\
\noindent 
E-mail: {\tt paolo.salani@unifi.it}\\

\smallskip
\noindent 
A. T.: Department of Mathematical Sciences,  Tokyo Metropolitan University, 
Minami-osawa, Hachioji-shi, Tokyo 192--0397, Japan\\
\noindent 
E-mail: {\tt asuka@tmu.ac.jp}\\
\newpage
\section{Introduction}
A  nonnegative function $u$ in ${\bf R}^N$ is said {\em logconcave} in ${\bf R}^N$ if 
$$
u((1-\mu)x+\mu x)\ge u(x)^{1-\mu}u(y)^\mu
$$
for $\mu\in[0,1]$ and $x$, $y\in{\bf R}^N$ such that $u(x)u(y)>0$.
This is equivalent to that 
the set $S_u:=\{x\in{\bf R}^N\,:\,u(x)>0\}$ is convex and
$\log u$ is concave in $S_u$. 
Logconcavity is a very useful variation of concavity and plays an important role 
in various fields such as PDEs, geometry, probability, statics, optimization theory and so on (see e.g. \cite{SW}). 
Most of its relevance, especially for elliptic and parabolic equations, is due to the fact that  
the Gauss kernel 
\begin{equation}
\label{eq:1.1}
G(x,t):=(4\pi t)^{-\frac{N}{2}}\exp\left(-\frac{|x|^2}{4t}\right)
\end{equation}
is logconcave in ${\bf R}^N$ for any fixed $t>0$. Indeed, 
\begin{equation}
\label{eq:1.2}
\log\,G(x,t)=-\frac{|x|^2}{4t}+\log(4\pi t)^{-\frac{N}{2}}
\end{equation}
is concave in ${\bf R}^N$ for any fixed $t>0$. 
Exploiting the logconcavity of the Gauss kernel, 
Brascamp and Lieb~\cite{BL} proved that logconcavity is preserved by the heat flow and they also obtained
the logconcavity of the first positive Dirichlet eigenfunction for the Laplace operator $-\Delta$ in a bounded convex domain. 
(See also \cite{DD, Kor}.) 
For later convenience, we state explicitly these two classical results below. 
\begin{itemize}
  \item[(a)] 
  Let $u$ be a bounded nonnegative solution of 
  \begin{equation}
  \label{eq:1.3}
  \left\{
  \begin{array}{ll}
  \partial_t u=\Delta u & \quad\mbox{in}\quad\Omega\times(0,\infty),\vspace{3pt}\\
  u=0 & \quad\mbox{on}\quad\partial\Omega\times(0,\infty)\quad\mbox{if}\quad\partial\Omega\not=\emptyset,\vspace{3pt}\\
  u(x,0)=u_0(x) & \quad\mbox{in}\quad\Omega,
  \end{array}
  \right.
  \end{equation}
  where $\Omega$ is a convex domain in ${\bf R}^N$ 
  and $u_0$ is a bounded nonnegative function in $\Omega$.
  Then $u(\cdot,t)$ is logconcave in $\Omega$ for any $t>0$ 
  if $u_0$ is logconcave in $\Omega$. 
  \item[(b)] 
  Let $\Omega$ be a bounded convex domain in ${\bf R}^N$ and $\lambda_1$ 
  the first Dirichlet eigenvalue for the Laplace operator in $\Omega$.
  If $\phi$ solves
  \begin{equation}
  \label{eq:1.4}
  \left\{\begin{array}{ll}
  -\Delta\phi=\lambda_1\,\phi\quad&\mbox{in}\quad\Omega,\vspace{3pt}\\
  \phi=0\quad&\mbox{on}\quad\partial\Omega,\vspace{3pt}\\
  \phi>0\quad&\mbox{in}\quad\Omega,
  \end{array}\right.
  \end{equation}
  then $\phi$ is logconcave in $\Omega$.
\end{itemize}
We denote by $e^{t\Delta_\Omega}u_0$ the (unique) solution to problem~\eqref{eq:1.3}. 
In particular, in the case of $\Omega={\bf R}^N$, 
we write $e^{t\Delta}u_0=e^{t\Delta_{{\bf R}^N}}u_0$, that is,
\begin{equation}
\label{eq:1.5}
[e^{t\Delta}u_0](x)=\int_{{\bf R}^N}G(x-y,t)u_0(y)\,dy,\quad x\in{\bf R}^N,\,\,t>0. 
\end{equation} 
Logconcavity is so naturally and deeply linked to heat transfer that 
$e^{t\Delta}u_0$ spontaneously becomes logconcave in ${\bf R}^N$ 
even without the logconcavity of initial function $u_0$. 
Indeed, Lee and V\'azquez~\cite{LV} proved the following: 
\begin{itemize}
  \item[(c)] 
  Let $u_0$ be a bounded nonnegative function in ${\bf R}^N$ with compact support. 
  Then there exists $T>0$ such that 
  $e^{t\Delta}u_0$ is logconcave in ${\bf R}^N$ for $t\ge T$.  (See \cite[Theorem~5.1]{LV}.)
\end{itemize}
Due to the above reasons, 
logconcavity is commonly regarded as the optimal concavity 
for the heat flow and for the first positive Dirichlet eigenfunction for $-\Delta$.
\vspace{7pt}

In this paper we dare to ask the following question:
\begin{itemize}
 \item[({\bf Q1})] 
 Is logconcavity the strongest concavity preserved by the heat flow in a convex domain? 
 If not, what is the strongest concavity preserved by the heat flow?
\end{itemize}
We introduce a new variation of concavity 
and give answers to ({\bf Q1}). 
More precisely, 
we introduce $\alpha$-logconcavity as a refinement of $p$-concavity at $p=0$ (see Section~2)
and show that the heat flow preserves 2-logconcavity (see Theorem~\ref{Theorem:3.1}). 
Here 2-logconcavity is stronger than usual logconcavity and 
we prove that 2-logconcavity is exactly the strongest concavity property preserved by the heat flow (see Theorem~\ref{Theorem:3.3}).
\medskip

Another natural question which spontaneously arises after ({\bf Q1}) is the following:
{\em is logconcavity the strongest concavity shared by the solution $\phi$ of \eqref{eq:1.4} for any bounded convex domain $\Omega$?}
We are not able to give here an exhaustive answer to this question, but we conjecture that it is negative and that also the first positive Dirichlet eigenfunction for $-\Delta$ in every convex domain is 2-logconcave. See Remark \ref{Theorem:3.2} about this. 
\medskip

The rest of this paper is organized as follows. 
In Section~2 we introduce a new variation of concavity and prove some lemmas. 
In particular, we show that 2-logconcavity is the strongest concavity for the Gauss kernel $G(\cdot,t)$ to satisfy. 
In Section~3 we state the main results of this paper. 
The proofs of the main results are given in Sections~4 and 5. 
\section{Logarithmic power concavity}
For $x\in{\bf R}^N$ and $R>0$, set $B(x,R):=\{y\in{\bf R}^N\,:\,|x-y|<R\}$. 
For any measurable set $E$, we denote by $\chi_E$ the characteristic function of $E$. 
Furthermore, for any function $u$ in a set $\Omega$ in ${\bf R}^N$, we say that $U$ is the {\it zero extension} of $u$ 
if $U(x)=u(x)$ for $x\in\Omega$ and $U(x)=0$ for $x\not\in\Omega$. 
We often identify $u$ with its zero extension~$U$. 
A function $f:{\bf R}^N\to{\bf R}\cup\{-\infty\}$ is a (proper) {\it concave} function if  $f((1-\mu)x+\mu y)\geq(1-\mu)f(x)+\mu f(y)$ for $x,y\in{\bf R}^N$, $\mu\in[0,1]$ (and $f(x)>-\infty$ for at least one $x\in{\bf R}^N$). Here we deal with $-\infty$ in the obvious way, that is: $-\infty+a=-\infty$ for $a\in{\bf R}$ and $-\infty\geq-\infty$.
\vspace{3pt}

Let $G=G(x,t)$ be the Gauss kernel (see \eqref{eq:1.1}). 
Similarly to \eqref{eq:1.2}, for any fixed $t>0$, 
it follows that 
$$
-\left(-\log\left[\kappa G(x,t)\right]\right)^{\frac{1}{2}}=-\left[\frac{|x|^2}{4t}-\log\left((4\pi t)^{-\frac{N}{2}}\kappa\right)\right]^{\frac{1}{2}}
$$
is still concave in ${\bf R}^N$ for any sufficiently small $\kappa>0$.
Motivated by this, we formulate a definition of $\alpha$-logconcavity $(\alpha>0)$. 
Let $L_\alpha=L_\alpha(s)$ be a strictly increasing function on $[0,1]$ 
defined by
$$
L_\alpha(s):=-(-\log s)^{\frac{1}{\alpha}}\quad\mbox{for}\quad s\in(0,1],\qquad
L_\alpha(s):=-\infty\quad\mbox{for}\quad s=0. 
$$
\begin{definition}
\label{Definition:2.1}
Let $\alpha>0$.
\vspace{3pt}
\newline
{\rm (i)} Let $u$ be a bounded nonnegative function in ${\bf R}^N$. 
We say that $u$ is $\alpha$-logconcave in ${\bf R}^N$ if 
$$
L_\alpha(\kappa u((1-\mu)x+\mu y))\geq (1-\mu)L_\alpha(\kappa u(x))+\mu L_\alpha(\kappa u(y)),
\quad x,y\in{\bf R}^N,\,\,\mu\in[0,1],
$$
for all sufficiently small $\kappa>0$. 
\vspace{3pt}
\newline
{\rm (ii)} Let $u$ be a bounded nonnegative function in a convex set $\Omega$ in ${\bf R}^N$. 
We say that $u$ is $\alpha$-logconcave in $\Omega$ if the zero extension of $u$ is $\alpha$-logconcave in ${\bf R}^N$. 
\end{definition}
Definition~\ref{Definition:2.1} means that, 
for any bounded nonnegative function $u$ in a convex set $\Omega$, 
$u$ is $\alpha$-logconcave in $\Omega$ if and only if 
the function~$L_\alpha(\kappa U)$ is concave in ${\bf R}^N$ for all sufficiently small $\kappa>0$. 
Due to Definition~\ref{Definition:2.1}, 
we easily see the following properties: 
\begin{itemize}
  \item 
  Logconcavity corresponds to $1$-logconcavity; 
  \item
  If $u$ is $\alpha$-logconcave in $\Omega$ for some $\alpha>0$, then 
  $\kappa u$ is also $\alpha$-logconcave in $\Omega$ for any $\kappa>0$;
  \item 
  If $0<\alpha\le\beta$ and $u$ is $\beta$-logconcave in $\Omega$, then $u$ is $\alpha$-logconcave in $\Omega$;
  \item
  The Gauss kernel $G(\cdot,t)$ is 2-logconcave in ${\bf R}^N$ for any $t>0$.
\end{itemize}
Furthermore, we have:
\begin{lemma}
\label{Lemma:2.1} 
Let $\alpha\ge 1$. 
Let $u$ be a function in a convex set $\Omega$ such that $0\le u\le 1$ in $\Omega$.  
If $L_\alpha(u)$ is concave in $\Omega$, 
then $L_\alpha(\kappa u)$ is also concave in $\Omega$ for any $0<\kappa\le 1$.
\end{lemma}
{\bf Proof.}
Let $x$, $y\in\Omega$ and $\mu\in[0,1]$. 
Assume that $u(x)u(y)>0$. 
Since $L_\alpha^{-1}(s)=\exp(-(-s)^\alpha)$ for $s\in(-\infty,0]$, 
we find
\begin{equation*}
\begin{split}
\psi(\kappa):=\, & \kappa^{-1}L_\alpha^{-1}
\left[(1-\mu)L_\alpha(\kappa u(x))+\mu L_\alpha(\kappa u(y))\right]\\
=\, & \kappa^{-1}\exp\left\{-\left[(1-\mu)(-\log \kappa u(x))^{\frac{1}{\alpha}}+\mu(-\log \kappa u(y))^{\frac{1}{\alpha}}\right]^\alpha\right\}
\end{split}
\end{equation*}
for $0<\kappa\le 1$. 
Since $L_\alpha(u)$ is concave in $\Omega$, it follows that
\begin{equation}
\label{eq:2.1}
(1-\mu)u(x)+\mu u(y)\ge\psi(1).
\end{equation}
For $a$, $b>0$ and $\gamma\in(-\infty,\infty)$, 
set 
$$
M_\gamma(a,b;\mu):=\left[(1-\mu)a^\gamma+\mu b^\gamma\right]^{\frac{1}{\gamma}}. 
$$
Then 
\begin{equation}
\label{eq:2.2}
\begin{split}
 & \psi'(\kappa)=-\kappa^{-2}\exp\left\{-M_{\frac{1}{\alpha}}(-\log \kappa u(x),-\log\kappa u(y);\mu)\right\}\\
 & \qquad\quad
 +\kappa^{-1}\exp\left\{-M_{\frac{1}{\alpha}}(-\log \kappa u(x),-\log\kappa u(y);\mu)\right\}\\
 & \qquad\qquad
 \times\left[(1-\mu)(-\log \kappa u(x))^{\frac{1}{\alpha}}+\mu(-\log \kappa u(y))^{\frac{1}{\alpha}}\right]^{\alpha-1}\\
 & \qquad\qquad
 \times\kappa^{-1}\left[(1-\mu)(-\log\kappa u(x))^{\frac{1-\alpha}{\alpha}}+\mu(-\log\kappa u(y))^{\frac{1-\alpha}{\alpha}}\right]\\
  & =\kappa^{-2}\exp\left\{-M_{\frac{1}{\alpha}}(-\log \kappa u(x),-\log\kappa u(y);\mu)\right\}\\
  & \times\left[-1+M_{\frac{1}{\alpha}}(-\log \kappa u(x),-\log\kappa u(y);\mu)^{\frac{\alpha-1}{\alpha}}
  M_{\frac{1-\alpha}{\alpha}}(-\log \kappa u(x),-\log\kappa u(y);\mu)^{\frac{1-\alpha}{\alpha}}\right]
\end{split}
\end{equation}
for $0<\kappa\le 1$. 
On the other hand, since $\alpha\ge 1$, 
it follows that $1/\alpha\ge(1-\alpha)/\alpha$. 
Then the Jensen inequality yields 
$$
M_{\frac{1}{\alpha}}(a,b;\mu)\ge M_{\frac{1-\alpha}{\alpha}}(a,b;\mu)\quad\mbox{for}\quad a,b>0. 
$$
This together with \eqref{eq:2.2} implies that $\psi'(\kappa)\ge 0$ for $0<\kappa\le 1$. 
Therefore, by \eqref{eq:2.1} we obtain 
\begin{equation}
\label{eq:2.3}
(1-\mu)u(x)+\mu u(y)\ge\psi(\kappa)
=\kappa^{-1}L_\alpha^{-1}
\left[(1-\mu)L_\alpha(\kappa u(x))+\mu L_\alpha(\kappa u(y))\right]
\end{equation}
for $0<\kappa\le 1$ in the case of $u(x)u(y)>0$. 
In the case of $u(x)u(y)=0$, 
by the definition of $L_\alpha$ we easily obtain \eqref{eq:2.3} for $0<\kappa\le 1$. 
These mean that $L_\alpha(\kappa u)$ is concave in $\Omega$ for $0<\kappa\le 1$. 
Thus Lemma~\ref{Lemma:2.1} follows. 
$\Box$
\begin{remark}
\label{Remark:2.1}
{\rm (i)} Let $u$ be a bounded nonnegative function in a convex set $\Omega$ and $\alpha>0$. 
We say that $u$ is weakly $\alpha$-logconcave in $\Omega$ if
$$
L_\alpha(\kappa U((1-\mu)x+\mu y))\geq (1-\mu)L_\alpha(\kappa U(x))+\mu L_\alpha(\kappa U(y)),
\quad x,y\in{\bf R}^N,\,\,\mu\in[0,1],
$$
for some $\kappa>0$. Here $U$ is the zero extension of $u$. 
Lemma~{\rm\ref{Lemma:2.1}} implies that 
$\alpha$-logconcavity is equivalent to weak $\alpha$-logconcavity in the case of $\alpha\ge 1$. 
\vspace{3pt}
\newline
{\rm (ii)} For $0<\alpha<1$, 
$\alpha$-logconcavity is not equivalent to weak $\alpha$-logconcavity. 
Indeed, set $u(x):=\exp(-|x|^\alpha)\chi_{B(0,R)}$, where $0<R\le\infty$. 
Then $L_\alpha(u)=-|x|$ is concave in $B(0,R)$. 
On the other hand, for any $0<\kappa<1$, we have 
\begin{equation*}
\begin{split}
\frac{\partial}{\partial r}L_\alpha (\kappa u(x)) & =-(-\log\kappa+|x|^\alpha)^{-1+\frac{1}{\alpha}}|x|^{-1+\alpha},\\
\frac{\partial^2}{\partial r^2}L_\alpha (\kappa u(x))
 & =(\alpha-1)(-\log\kappa+|x|^\alpha)^{-2+\frac{1}{\alpha}}|x|^{-2+2\alpha}
+(1-\alpha)(-\log\kappa+|x|^\alpha)^{-1+\frac{1}{\alpha}}|x|^{-2+\alpha}\\
 & =(1-\alpha)(-\log\kappa+|x|^\alpha)^{-2+\frac{1}{\alpha}}|x|^{-2+\alpha}(-\log\kappa)>0,
\end{split}
\end{equation*}
for $x\in B(0,R)\setminus\{0\}$, where $r:=|x|>0$. This means that $L_\alpha (\kappa u)$ is not concave in $B(0,R)$ for any $0<\kappa<1$. 
\end{remark}

Next we introduce the notion of $F$-concavity, which generalises and embraces all the notions of concavity we have already seen.
\begin{definition}
\label{Definition:2.2}
Let $\Omega$ be a convex set in ${\bf R}^N$.
\newline
{\rm (i)} A function $F:[0,1]\to{\bf R}\cup\{-\infty\}$ is said admissible 
if $F$ is strictly increasing continuous in $(0,1]$, $F(0)=-\infty$ and $F(s)\not=-\infty$ for $s>0$. 
\vspace{3pt}
\newline
{\rm (ii)} 
Let $F$ be admissible. 
Let $u$ be a bounded nonnegative function in $\Omega$ 
and $U$ the zero extension of $u$. 
Then $u$ is said $F$-concave in $\Omega$ if 
$0\le\kappa U(x)\le 1$ in ${\bf R}^N$ and 
$$
F(\kappa U((1-\mu)x+\mu y))\geq (1-\mu)F(\kappa U(x))+\mu F(\kappa U(y)),
\quad x,y\in{\bf R}^N,\,\,\mu\in[0,1],
$$ 
for all sufficiently small $\kappa>0$. 
We denote by ${\mathcal C}_\Omega[F]$ 
the set of $F$-concave functions in $\Omega$. 
Furthermore, in the case of $\Omega={\bf R}^N$, 
we write ${\mathcal C}[F]={\mathcal C}_\Omega[F]$ for simplicity. 
\vspace{3pt}
\newline
{\rm (iii)} 
Let $F_1$ and $F_2$ be admissible. 
We say that
$F_1$-concavity is stronger than $F_2$-concavity in $\Omega$ 
if ${\mathcal C}_\Omega[F_1]\subsetneq{\mathcal C}_\Omega[F_2]$. 
\end{definition}
We recall that 
a bounded nonnegative function $u$ in a convex set $\Omega$ is said {\em $p$-concave} in $\Omega$, where $p\in{\bf R}$, 
if $u$ is $F$-concave with $F=F_p$ in $\Omega$, where  
$$
F_p(s):=
\left\{
\begin{array}{ll}
\displaystyle{\frac{1}{p}s^p} & \mbox{for $s>0$ if $p\not=0$},\vspace{7pt}\\
\log s & \mbox{for $s>0$ if $p=0$},\vspace{5pt}\\
-\infty & \mbox{for $s=0$}.
\end{array}
\right.
$$
Here $1$-concavity corresponds to usual concavity while $0$-concavity corresponds to usual logconcavity (in other words, $1$-logconcavity).
Furthermore, 
$u$ is said {\em quasiconcave} or {\em $-\infty$-concave} in $\Omega$ 
if all superlevel sets of $u$ are convex, 
while it is said {\em $\infty$-concave} in $\Omega$ 
if $u$ satisfies 
$$
u((1-\mu)x+\mu y)\ge\max\{u(x),u(y)\}
$$
for $x$, $y\in\Omega$ with $u(x)u(y)>0$ and $\mu\in[0,1]$.
Then, by the Jensen inequality we have:
\begin{itemize}
  \item Let $-\infty\le p\le q\le\infty$. 
  If $u$ is $q$-concave in a convex set $\Omega$, then $u$ is also $p$-concave in~$\Omega$.
  \end{itemize}
Among concavity properties, apart from usual concavity, of course logconcavity has been the most deeply investigated, especially for its importance in probability and convex geometry (see for instance \cite{Co} for an overview and the series of papers \cite{AKM, AM1, AM2,M},
which recently broadened and structured the theory of log-concave functions). 
Clearly, if a function $u$ is $F$-concave in $\Omega$ for some admissible $F$, then it is quasiconcave in $\Omega$; 
vice versa, if $u$ is $\infty$-concave in $\Omega$, then it is $F$-concave in $\Omega$ for any admissible $F$. 
These mean that quasiconcavity (resp.\,\,$\infty$-concavity) is  
the weakest (resp.\,\,strongest) conceivable concavity.  
Notice that $\alpha$-logconcavity ($\alpha>0$) corresponds to $F$-concavity with $F=L_\alpha$ 
and it is weaker (resp.\,\,stronger) than $p$-concavity for any $p>0$ (resp.\,\,$p<0$). 
Indeed, the following lemma holds. 
\begin{lemma}
\label{Lemma:2.2}
Let $\Omega$ be a convex set in ${\bf R}^N$ 
and $u$ a nonnegative bounded function in $\Omega$.
\vspace{3pt}
\newline
{\rm (i)} 
If $u$ is $p$-concave in $\Omega$ for some $p>0$, then 
$u$ is $\alpha$-logconcave in $\Omega$ for any $\alpha>0$. 
\vspace{3pt}
\newline
{\rm (ii)} 
If $u$ is $\alpha$-logconcave in $\Omega$ for some $\alpha>0$, 
then $u$ is $p$-concave in $\Omega$ for any $p<0$. 
\end{lemma}
{\bf Proof.}
We prove assertion~(i). 
Let $p>0$ and $\alpha>0$. It suffices to prove that  
\begin{equation}
\label{eq:2.4}
\left[(1-\mu)a^p+\mu b^p\right]^{\frac{1}{p}}
\ge\exp\left\{-\left[(1-\mu)(-\log a)^{\frac{1}{\alpha}}+\mu(-\log b)^{\frac{1}{\alpha}}\right]^\alpha\right\}
\end{equation}
holds for all sufficiently small $a$, $b>0$ and all $\mu\in[0,1]$. 
This is equivalent to that the inequality 
\begin{equation}
\label{eq:2.5}
\left(-\frac{1}{p}\log\left[(1-\mu)\tilde{a}+\mu\tilde{b}\right]\right)^{\frac{1}{\alpha}}
\le(1-\mu)\left(-\frac{1}{p}\log\tilde{a}\right)^{\frac{1}{\alpha}}+\mu\left(-\frac{1}{p}\log\tilde{b}\right)^{\frac{1}{\alpha}}
\end{equation}
holds for all sufficiently small $\tilde{a}:=a^p$, $\tilde{b}:=b^p>0$ and all $\mu\in[0,1]$. 
Inequality~\eqref{eq:2.5} follows from the fact that 
the function 
$$
s\mapsto\left(-\frac{1}{p}\log s\right)^{\frac{1}{\alpha}} 
$$
is convex for all sufficiently small $s>0$. 
Thus \eqref{eq:2.4} holds for all sufficiently small $a$, $b>0$ and all $\mu\in[0,1]$ 
and assertion~(i) follows. 
Similarly, we obtain assertion~(ii) and the proof is complete. 
$\Box$\vspace{5pt}
\newline
Lemma~\ref{Lemma:2.2} implies that 
$\alpha$-logconcavity is a refinement of $p$-concavity at $p=0$.
\vspace{5pt}

At the end of this section 
we show that 2-logconcavity is the strongest concavity for the Gauss kernel $G(\cdot,t)$ to satisfy. 
This plays a crucial role in giving an answer to the second part of ({\bf Q1}). 
\begin{lemma}
\label{Lemma:2.3}
Let $F$ be admissible such that 
$G(\cdot,t)$ is $F$-concave in ${\bf R}^N$ for some $t>0$. 
Then a bounded nonnegative function $u$ in ${\bf R}^N$ is $F$-concave in ${\bf R}^N$ if $u$ is $2$-logconcave in ${\bf R}^N$ $($in other words 
${\mathcal C}[L_2]\subset{\mathcal C}[F]$$)$. 
Furthermore, 
\begin{equation}
\label{eq:2.6}
{\mathcal C}[L_2]
=\bigcap_{F\in\{H\,:\,G(\cdot,t)\in{\mathcal C}[H]\}}{\mathcal C}[F]
\quad
\mbox{for any $t>0$}.
\end{equation}
\end{lemma}
{\bf Proof.} 
Assume that $G(\cdot,t)$ is $F$-concave in ${\bf R}^N$ for some $t>0$. 
It follows from Definition~\ref{Definition:2.1} that 
the function $e^{-|x|^2}$ is $F$-concave in ${\bf R}^N$. 
Then we obtain the $F$-concavity of $e^{-s^2}$ $(s\in{\bf R})$. 

Let $u$ be 2-logconcave in ${\bf R}^N$. By Definition~\ref{Definition:2.1} 
we see that $L_2(\kappa u)$ is concave in ${\bf R}^N$ for all sufficiently small $\kappa>0$. 
Set 
$$
w(x):=-L_2(\kappa u(x))
=\left\{
\begin{array}{ll}
\sqrt{-\log\kappa u(x)}\quad & \mbox{if}\quad u(x)>0,\vspace{3pt}\\
\infty\quad & \mbox{if}\quad u(x)=0.
\end{array}
\right.
$$
Then $w$ is nonnegative and convex in ${\bf R}^N$, that is,
\begin{equation}
\label{eq:2.7}
0\le w((1-\mu)x+\mu y)\le(1-\mu)w(x)+\mu w(y)
\end{equation}
for $x$, $y\in{\bf R}^N$ and $\mu\in[0,1]$.  
On the other hand, 
by $F$-concavity of $e^{-s^2}$  
we have
\begin{equation}
\label{eq:2.8}
 F\left(\kappa e^{-[(1-\mu)w(x)+\mu w(y)]^2}\right)
 \ge(1-\mu)F(\kappa e^{-w(x)^2})+\mu F(\kappa e^{-w(y)^2})
\end{equation}
for all sufficiently small $\kappa>0$. 
Since $F$ is an increasing function, 
by \eqref{eq:2.7} and \eqref{eq:2.8} we obtain 
\begin{equation*}
\begin{split}
 & F(\kappa^2 u((1-\mu)x+\mu y))=F\left(\kappa\exp(-w((1-\mu)x+\mu y)^2)\right)\\
 & \qquad
 \ge F\left(\kappa e^{-[(1-\mu)w(x)+\mu w(y)]^2}\right)
\ge(1-\mu)F(\kappa e^{-w(x)^2})+\mu F(\kappa e^{-w(y)^2})\\
 & \qquad
 =(1-\mu)F(\kappa^2 u(x))+\mu F(\kappa^2 u(y))
\end{split}
\end{equation*}
for all sufficiently small $\kappa>0$ 
if $u(x)u(y)>0$. 
This inequality also holds in the case of $u(x)u(y)=0$. 
These imply that $u$ is $F$-concave in ${\bf R}^N$ and  
\begin{equation}
\label{eq:2.9}
{\mathcal C}[L_2]
\subset\,\bigcap_{F\in\{H\,:\,G(\cdot,t)\in{\mathcal C}[H]\}}
{\mathcal C}[F]. 
\end{equation}
On the other hand, since $G(\cdot,t)$ is 2-logconcave, 
it turns out that 
$$
\bigcap_{F\in\{H\,:\,G(\cdot,t)\in{\mathcal C}[H]\}}
{\mathcal C}[F]\subset {\mathcal C}[L_2].
$$
This together with \eqref{eq:2.9} implies \eqref{eq:2.6}. 
Thus Lemma~\ref{Lemma:2.3} follows. 
$\Box$
\section{Main results}
We are now ready to state the main results of this paper.
The first one ensures that 
the heat flow preserves $\alpha$-logconcavity with $1\le\alpha\le 2$. 
\begin{theorem}
\label{Theorem:3.1}
Let $\Omega$ be a convex domain in ${\bf R}^N$ and $1\le\alpha\le 2$. 
Let $u_0$ be a bounded nonnegative function in $\Omega$ and $u:=e^{t\Delta_\Omega}u_0$. 
Assume that $0\le u_0\le 1$ and $L_\alpha(u_0)$ is concave in $\Omega$. 
Then $L_\alpha(u(\cdot,t))$ is concave in $\Omega$ for any $t>0$.
\end{theorem}
Since $\alpha$-logconcavity with $\alpha>1$ is stronger than usual logconcavity, 
Theorems~\ref{Theorem:3.1} gives answer to 
the first part of ({\bf Q1}).
Furthermore, as a corollary of Theorem~\ref{Theorem:3.1}, we have the following.
\begin{corollary}
\label{Corollary:3.1}
Let $\Omega$ be a convex domain in ${\bf R}^N$. 
Let $u_0$ be a bounded nonnegative function in $\Omega$. 
If $1\le\alpha\le 2$ and $u_0$ is $\alpha$-logconcave in $\Omega$, 
then $e^{t\Delta_\Omega}u_0$ is $\alpha$-logconcave in $\Omega$ for any $t>0$. 
\end{corollary}
Next we state a result which shows that 
$2$-logconcavity is the strongest concavity preserved by the heat flow. 
This addresses the second part of ({\bf Q1}).
\begin{theorem}
\label{Theorem:3.3} 
Let $F$ be admissible 
and $\Omega$ a convex domain in ${\bf R}^N$. 
Assume that $F$-concavity is stronger than {\rm 2}-logconcavity in $\Omega$, that is,  
${\mathcal C}_\Omega[F]\subsetneq{\mathcal C}_\Omega[L_2]$. 
Then there exists $u_0\in{\mathcal C}_\Omega[F]$ such that 
$$
e^{T\Delta_\Omega}u_0\notin{\mathcal C}_\Omega[F]\quad\mbox{for some $T>0$.}
$$
\end{theorem}
Here the following question naturally arises: 
\begin{itemize}
  \item[({\bf Q2})] 
 What is the weakest concavity preserved by the heat flow? 
\end{itemize}
Unfortunately we have no answers to ({\bf Q2}) and it is open. 
Notice that the heat flow does not necessarily preserve $p$-concavity for some $p<0$. 
See \cite{IS1, IS2}. (See also \cite{CW}.)  
%

Finally we assure that 
$e^{t\Delta}u_0$ spontaneously becomes $\alpha$-logconcave  
for any $\alpha\in[1,2)$ if $u_0$ has compact support. 
This improves assertion~(c).
\begin{theorem}
\label{Theorem:3.4} 
Let $u_0$ be a bounded nonnegative function in ${\bf R}^N$ 
with compact support. 
Then, for any given $1\le \alpha<2$, 
there exists $T_\alpha>0$ such that, 
for any $t\ge T_\alpha$, 
$L_\alpha(e^{t\Delta}u_0)$ is concave in ${\bf R}^N$, 
in particular, $e^{t\Delta}u_0$ is $\alpha$-logconcave in ${\bf R}^N$. 
\end{theorem}
We conjecture that Theorem~{\rm\ref{Theorem:3.4}} holds true even for $\alpha=2$, but 
we can not prove it here. Indeed, in our proof of Theorem~{\rm\ref{Theorem:3.4}}, 
$T_\alpha\to\infty$ as $\alpha\to 2$. 

\vspace{5pt}

In Section~4 we prove Theorems~\ref{Theorem:3.1} and \ref{Theorem:3.3}. 
Theorem~\ref{Theorem:3.1} is shown as an application of \cite{INS} 
however the proof is somewhat tricky (see Remark~\ref{Remark:4.1}). 
Furthermore, we prove Theorem~\ref{Theorem:3.3} by the use of Lemma~\ref{Lemma:2.3}.
In Section 5 we study the large time behavior of the second order derivatives of $e^{t\Delta}u_0$.  
This proves Theorem~\ref{Theorem:3.4}. 
\section{Proofs of Theorems~\ref{Theorem:3.1} and \ref{Theorem:3.3}}
Firstly we prove Theorem~\ref{Theorem:3.1} 
and show the preservation of $\alpha$-logconcavity $(1\le\alpha\le 2)$ by the heat flow. 
\vspace{3pt}
\newline
{\bf Proof of Theorem~\ref{Theorem:3.1}.}
Let $\Omega$ be a convex domain in ${\bf R}^N$. 
Let $u_0$ be a nontrivial function in $\Omega$ such that $0\le u_0(x)\le 1$ in $\Omega$. 
Then it follows from the strong maximum principle that $0<u<1$ in $\Omega\times(0,\infty)$. 
Assume that $L_\alpha(u_0)$ is concave in $\Omega$ for some $\alpha\in[1,2]$. 
\vspace{3pt}
\newline
\underline{1st step}: 
We consider the case where $\Omega$ is a bounded smooth convex domain, 
$u_0\in C(\overline\Omega)$ and $u_0=0$ on $\partial\Omega$. 
Set 
$$
w(x,t):=-L_\alpha(u(x,t))\ge 0,\qquad w_0(x):=-L_\alpha(u_0(x))\ge 0.
$$
Here $w_0$ is convex in $\Omega$. 
Then it follows that 
\begin{equation}
\label{eq:4.1}
\left\{\begin{array}{ll}
w_t-\Delta w+\displaystyle{\frac{1}{\gamma}\frac{|\nabla w|^2}{w^{\frac{\gamma-1}{\gamma}}}}
\displaystyle{+\frac{\gamma-1}{\gamma}\frac{|\nabla w|^2}{w}}=0 & \mbox{in}\quad\Omega\times(0,\infty)\,,\vspace{5pt}\\
w(x,0)=w_0(x) & \mbox{in}\quad\Omega\,,\vspace{7pt}\\
w>0 & \mbox{in}\quad\Omega\times(0,\infty)\,,\vspace{7pt}\\
w(x,t)\to+\infty & \mbox{as}\quad\mbox{$\mbox{dist}(x,\partial\Omega)\to 0$ for any $t>0$}\,,
\end{array}\right.
\end{equation}
where $\gamma:=1/\alpha\in[1/2,1]$. 

We prove that $w(\cdot,t)$ is convex in $\Omega$ for any $t>0$. 
For this aim, we set $z:=e^{-w}$ and show that $z(\cdot,t)$ is logconcave in $\Omega$ for any $t>0$. 
(See Remark~\ref{Remark:4.1}.) 
It follows from \eqref{eq:4.1} that
\begin{equation*}
\left\{\begin{array}{ll}
z_t-\Delta z+\displaystyle{\frac{|\nabla z|^2}{z}}
\left[-\frac{1}{\gamma}(-\log z)^{-\frac{\gamma-1}{\gamma}}+\frac{\gamma-1}{\gamma}(\log z)^{-1}+1\right]=0
 & \mbox{in}\quad\Omega\times(0,\infty)\,,\\
\\
z(x,0)=e^{-w_0(x)} & \mbox{in}\quad \Omega\,,\\
\\ 
z(x,t)=0\quad& \mbox{on}\quad\partial\Omega\times(0,\infty)\,.
\end{array}\right.
\end{equation*}
Furthermore, thanks to the convexity of $w_0$, we see that $z(\cdot,0)=e^{-w_0}$ is logconcave in $\Omega$. 
Applying \cite[Theorem 4.2, Corollary 4.2]{INS} (see also \cite[Theorem~4.2]{GK}),  
we deduce that 
$z(\cdot,t)$ is logconcave in $\Omega$ for any $t>0$ if
\begin{equation*}
\begin{split}
 & \mbox{$h(s,A):=e^{-s}\left[-e^s\,\text{trace}(A)+f(e^s,e^s\theta)\right]$ is convex}\\
 & \text{for $(s,A)\in(-\infty,0)\times {\rm Sym}_N$ for any fixed $\theta\in{\bf R}^N$}.
\end{split}
\end{equation*}
Here ${\rm Sym}_N$ denotes the space of real $N\times N$ symmetric matrices  
and 
$$
f(\zeta,\vartheta):=\frac{|\vartheta|^2}{\zeta}\left[-\frac{1}{\gamma}(-\log \zeta)^{-\frac{\gamma-1}{\gamma}}+\frac{\gamma-1}{\gamma}(\log \zeta)^{-1}+1\right]
\quad\mbox{for}\quad (\zeta,\vartheta)\in(0,1)\times{\bf R}^N.
$$
On the other hand, for any fixed $\theta\in{\bf R}^N$, 
$$
h(s,A)=-\text{trace}(A)+|\theta|^2
\left[-\frac{1}{\gamma}(-s)^{-\frac{\gamma-1}{\gamma}}+\frac{\gamma-1}{\gamma}s^{-1}+1\right]
$$
is convex for $(s,A)\in(-\infty,0)\times {\rm Sym}_N$ if and only if $1/2\le\gamma\le 1$. 
Therefore $z(\cdot,t)$ is logconcave in $\Omega$ for any $t>0$. 
This implies that 
$u(\cdot,t)$ is $\alpha$-logconcave for any $t>0$. 
\vspace{3pt}
\newline
\underline{2nd step}: 
We consider the case where $\Omega$ is a bounded smooth convex domain. 
In this step we do not assume that $u_0=0$ on $\partial\Omega$. 
Since $u_0$ is $\alpha$-logconcave in $\Omega$, 
we see that $L_\alpha(u_0)$ is concave in $\Omega$. 
Set 
$$
v_0(x):=
\left\{
\begin{array}{ll}
\exp(L_\alpha(u_0(x))) & \mbox{for}\quad x\in\Omega,\vspace{3pt}\\
0 & \mbox{for}\quad x\not\in\Omega,
\end{array}
\right.
\qquad
v(x,t):=[e^{t\Delta}v_0](x),
$$
for $x\in{\bf R}^N$ and $t>0$. 
Here we let $e^{-\infty}:=0$. 
Then $v_0$ is logconcave in ${\bf R}^N$. 
We deduce from assertion~(a) that $v(\cdot,t)$ is logconcave for any $t>0$. 
Furthermore, we deduce from $v_0\in L^1({\bf R}^N)\cap L^\infty({\bf R}^N)$ that 
\begin{eqnarray}
\label{eq:4.2}
 & & \mbox{$v(\cdot, t)$ is a positive continuous function in ${\bf R}^N$ for any $t>0$},\\  
\label{eq:4.3}
 & & \|v(t)\|_{L^\infty({\bf R}^N)}<\|v_0\|_{L^\infty({\bf R}^N)}\mbox{ for any $t>0$},\\
\label{eq:4.4}
 & & \lim_{t\to 0}\|v(t)-v_0\|_{L^1({\bf R}^N)}=0.
\end{eqnarray}
By \eqref{eq:4.4} we can find a sequence $\{t_n\}\subset(0,\infty)$ with $\lim_{n\to\infty}t_n=0$ such that 
\begin{equation}
\label{eq:4.5}
\lim_{n\to\infty}v(x,t_n)=v_0(x)
\end{equation}
for almost all $x\in{\bf R}^N$. 

Let $\eta$ solve
$$
-\Delta \eta=1\quad\mbox{in}\quad\Omega,\qquad \eta>0\quad\mbox{in}\quad\Omega,
\qquad \eta=0\quad\mbox{on}\quad\partial\Omega.
$$
Then $\eta$ is $1/2$-concave in $\Omega$ (see e.g. \cite[Theorem~4.1]{Ken}), which implies that 
$\log \eta$ is concave in $\Omega$ and $\log \eta\to-\infty$ as $\mbox{dist}\,(x,\partial\Omega)\to 0$. 
By \eqref{eq:4.2} and \eqref{eq:4.3} we can find a sequence $\{m_n\}\subset(1,\infty)$ 
with $\lim_{n\to\infty}m_n=\infty$ such that 
$$
V_n(x):=\log v(x,t_n)+m_n^{-1}\log \eta(x)
$$
is continuous and concave in $\Omega$ and 
$$
\sup_{x\in\Omega}V_n(x)\le
\underset{x\in\Omega}{\mbox{ess sup}}\,\log v_0.
$$ 
Furthermore, by \eqref{eq:4.5} we have 
\begin{equation*}
\begin{split}
 & \lim_{n\to\infty} V_n(x)=\log v_0(x)=L_\alpha(u_0(x))\quad\mbox{for almost all $x\in\Omega$},\\
 & V_n(x)\to-\infty\quad\mbox{as}\quad \mbox{dist}\,(x,\partial\Omega)\to 0.
\end{split}
\end{equation*}
Then the function 
$u_{0,n}(x):=L_\alpha^{-1}(V_n(x))$
satisfies 
\begin{equation}
\label{eq:4.6}
\mbox{$0\le u_{0,n}\le 1$ in $\Omega$, $u_{0,n}=0$ on $\partial\Omega$ and 
 $\displaystyle{\lim_{n\to\infty}}u_{0,n}(x)=u_0(x)$ for almost all $x\in\Omega$}.
\end{equation}
Furthermore, $u_{0,n}$ is continuous on $\overline{\Omega}$ 
and $L_\alpha(u_{0,n})$ is concave in $\Omega$. 
Let 
$$
u_n(x,t):=[e^{t\Delta_\Omega}u_{0,n}](x)=\int_\Omega G_\Omega (x,y,t)u_{0,n}(y)\,dy,
$$
where $G_\Omega=G_\Omega(x,y,t)$ is the Dirichlet heat kernel in $\Omega$. 
Then, by \eqref{eq:4.6} we apply the Lebesgue dominated convergence theorem to obtain
\begin{equation}
\label{eq:4.7}
\lim_{n\to\infty}u_n(x,t)=\int_\Omega G_\Omega (x,y,t)u_0(y)\,dy=u(x,t),
\quad x\in\Omega,\,\,t>0. 
\end{equation}
On the other hand, by the argument in 1st step we see that $L_\alpha(u_n(\cdot,t))$ is concave in $\Omega$ for any $t>0$. 
Then we deduce from \eqref{eq:4.7} that 
$L_\alpha(u(\cdot,t))$ is also concave in $\Omega$ for any $t>0$. 
Thus Theorem~\ref{Theorem:3.1} follows in the case where $\Omega$ is a bounded smooth convex domain.
\vspace{3pt}
\newline
\underline{3rd step}:
We complete the proof of Theorem~\ref{Theorem:3.1}.
There exists a sequence of bounded convex smooth domains $\{\Omega_n\}_{n=1}^\infty$ 
such that 
$$
\Omega_1\subset\Omega_2\subset\cdots\subset\Omega_n\subset\cdots,
\qquad
\bigcup_{n=1}^\infty\Omega_n=\Omega.
$$
(This is for instance a trivial consequence of \cite[Theorem 2.7.1]{sc}).

For any $n=1,2,\dots$, let $u_n:=e^{t\Delta_{\Omega_n}}(u_0\chi_{\Omega_n})$.  
The argument in 2nd step implies that 
$L_\alpha(u_n(\cdot,t))$ is concave in $\Omega_n$ for any $t>0$. 
Furthermore, 
by the comparison principle 
we see that 
\begin{equation*}
\begin{split}
 & u_n(x,t)\le u_{n+1}(x,t)\le u(x,t)\quad\mbox{in}\quad\Omega_n\times(0,\infty),\\
 & u(x,t)=\lim_{n\to\infty}u_n(x,t)\qquad\qquad\,\,\,\mbox{in}\quad\Omega\times(0,\infty).
\end{split}
\end{equation*}
Then we observe that $L_\alpha(u(\cdot,t))$ is concave in $\Omega$ for any $t>0$. 
Thus Theorem~\ref{Theorem:3.1} follows. 
$\Box$
\begin{remark} 
\label{Remark:4.1}
Sufficient conditions for the concavity of solutions to parabolic equations 
were discussed in {\rm\cite[{\it Section}~4.1]{INS}}. 
However we can not apply the arguments in {\rm\cite[{\it Section}~4.1]{INS}} 
to show the concavity of $-w(\cdot,t)$,
because assumption~{\rm (F3)} with $p=1$ in {\rm\cite{INS}} is not satisfied for the equation satisfied by  $-w$.  
\end{remark}
\begin{remark} 
\label{Theorem:3.2}
Theorem~{\rm\ref{Theorem:3.1}} implies that $e^{t\Delta_\Omega}\chi_\Omega(x)$ is $2$-logconcave (with respect to $x$) for every $t>0$.
As it is well known, 
by the eigenfunction expansion of solutions and the regularity theorems for the heat equation,
we have
$$
\lim_{t\to\infty}e^{\lambda_1 t}[e^{t\Delta_\Omega}\chi_\Omega](x)=c\phi(x)/\|\phi\|_{L^2(\Omega)}
$$
uniformly on $\overline{\Omega}$, 
where $\lambda_1$ and $\phi$ are as in assertion~{\rm (b)} of the Introduction and
$$
c=\int_\Omega \phi(x)\,dx\biggr/\|\phi\|_{L^2(\Omega)}>0\,.
$$
Then we may think to obtain the $2$-logcocanvity of $\phi$ just by letting $t\to+\infty$ and using the preservation of $2$-logconcavity by pointwise convergence. Unfortunately this approach does not work, since the parameter $\kappa$ of Definition~{\rm\ref{Definition:2.1}} for $e^{t\Delta_\Omega}\chi_\Omega(x)$ may tend to $0$ as 
$t$ tends to $+\infty$, while $2$-logconcavity is preserved only if $\kappa$ remains strictly positive.
\end{remark}
{\bf Proof of Corollary~\ref{Corollary:3.1}.}
Corollary~\ref{Corollary:3.1} directly follows from Theorem~\ref{Theorem:3.1}, 
Definition~\ref{Definition:2.1} and the linearity of the heat equation. 
$\Box$
\medskip

At the end of this section 
we prove Theorem~\ref{Theorem:3.3} with the aid of Lemma~\ref{Lemma:2.3}.
\vspace{5pt}
\newline
{\bf Proof of Theorem~\ref{Theorem:3.3}.}
Let us consider the case of $\Omega={\bf R}^N$. 
Since $F$ is stronger than $2$-logconcavity, 
by Lemma~\ref{Lemma:2.3} we see that $G(\cdot,t)\not\in{\mathcal C}[F]$ for any $t>0$. 
Then, for any $\epsilon>0$, 
there exist $\kappa\in(0,\epsilon)$, $\mu\in[0,1]$ and $x_1,x_2\in{\bf R}^N$ such that 
\begin{equation}
\label{eq:4.8}
F\left(\kappa e^{-|(1-\mu)x_1+\mu x_2|^2}\right)
-(1-\mu)F\left(\kappa e^{-|x_1|^2}\right)-\mu F\left(\kappa e^{-|x_2|^2}\right)<0\,.
\end{equation}
Let $K$ be a bounded convex set in ${\bf R}^N$ such that $|K|>0$ and set $u:=e^{t\Delta}\chi_K$. 
Since $\chi_K$ is $\infty$-concave, we see that $\chi_K$ is $F$-concave. 

On the other hand, 
it follows from \eqref{eq:1.5} that 
$$
\lim_{t\to\infty}t^{\frac{N}{2}}\|u(t)-|K|G(t)\|_{L^\infty({\bf R}^N)}=0. 
$$
This implies that 
\begin{equation}
\label{eq:4.9}
\lim_{t\to\infty}(4\pi t)^{\frac{N}{2}}|K|^{-1}u(2\sqrt{t}\xi,t)=e^{-|\xi|^2},\qquad \xi\in{\bf R}^N.
\end{equation}
For any $t\ge 1$ and $i=1,2$, set $\xi_i^t:=2\sqrt{t}x_i$. 
Since $F$ is continuous in $(0,1]$, by \eqref{eq:4.8} and \eqref{eq:4.9} we have  
\begin{equation*}
\begin{split}
 & F\left((4\pi t)^{\frac{N}{2}}|K|^{-1}\kappa u((1-\mu)\xi_1^t+\mu \xi_2^t,t)\right)\\
 & \qquad\quad
 -(1-\mu)F\left((4\pi t)^{\frac{N}{2}}|K|^{-1}\kappa u(\xi_1^t,t)\right)
 -\mu F\left((4\pi t)^{\frac{N}{2}}|K|^{-1}\kappa u(\xi_2^t,t)\right)\\
 & \to F\left(\kappa e^{-|(1-\mu)x_1+\mu x_2|^2}\right)-(1-\mu)F\left(\kappa e^{-|x_1|^2}\right)-\mu F\left(\kappa e^{-|x_2|^2}\right)<0
\end{split}
\end{equation*}
as $t\to\infty$.
Since $\epsilon$ is arbitrary, we see that 
\begin{equation}
\label{eq:4.10}
\mbox{$u(\cdot,T)=e^{T\Delta}\chi_K$ is not $F$-concave for all sufficiently large $T$.}
\end{equation}
Thus $F$-concavity is not preserved by the heat flow in ${\bf R}^N$. 

Next we consider the case of $\Omega\not={\bf R}^N$. 
We can assume, without loss of generality, that $0\in\Omega$ 
and $K\subset\Omega$. 
For $n=1,2,\dots$, set $\Omega_n:=n\Omega$. 
Then 
\begin{equation}
\label{eq:4.11}
\lim_{n\to\infty}[e^{t\Delta_{\Omega_n}}\chi_K](x)=[e^{t\Delta}\chi_K](x)
\end{equation}
for any $x\in{\bf R}^N$ and $t>0$. 
Let $T'$ be a sufficiently large constant. 
By \eqref{eq:4.10} and \eqref{eq:4.11} 
we observe that 
$e^{T'\Delta_{\Omega_n}}\chi_K$ is not $F$-concave for all sufficiently large $n$. 
Since 
$$
[e^{n^2t\Delta_{\Omega_n}}\chi_K](nx)=[e^{t\Delta_\Omega}\chi_{n^{-1}K}](x),
\qquad x\in\Omega,\,\,t>0,
$$
we see that 
$e^{n^{-2}T'\Delta_\Omega}\chi_{n^{-1}K}$ is not $F$-concave for all sufficiently large $n$. 
Combining the fact that $\chi_{n^{-1}K}$ is $F$-concave, we see that 
$F$-concavity is not preserved by the heat flow in $\Omega$. 
Thus Theorem~\ref{Theorem:3.3} follows. 
$\Box$
\section{Proof of Theorem~\ref{Theorem:3.4}}
We modify the arguments in the proof of \cite[Theorem~5.1]{LV} 
and prove Theorem~\ref{Theorem:3.4}.
\vspace{3pt}
\newline
{\bf Proof of Theorem~\ref{Theorem:3.4}.}
Let $u_0$ be a nontrivial bounded nonnegative function in ${\bf R}^N$ such that 
$\mbox{supp}\,u_0\subset B(0,R)$ for some $R>0$. 
Without loss of generality, we can assume that 
\begin{equation}
\label{eq:5.1}
\int_{{\bf R}^N}y_i u_0(y)\,dy=0,\qquad i=1,\dots,N. 
\end{equation}
Let $u:=e^{t\Delta}u_0$. It follows from \eqref{eq:1.1} and \eqref{eq:1.5} that 
\begin{equation}
\label{eq:5.2}
0<u\le\min\left\{\|u_0\|_{L^\infty({\bf R}^N)},(4\pi t)^{-\frac{N}{2}}\|u_0\|_{L^1({\bf R}^N)}\right\}
\end{equation}
for $(x,t)\in{\bf R}^N\times(0,\infty)$. In particular, 
$0<u(x,t)<1$ in ${\bf R}^N\times(T,\infty)$ for some $T>0$.

Let $\gamma:=1/\alpha\in(1/2,1]$. For the proof of Theorem~\ref{Theorem:3.4}, 
it suffices to prove that 
$$
v(x,t):=(-\log u(x,t))^\gamma
$$ 
is convex in ${\bf R}^N$ for all sufficiently large $t$.   
By \eqref{eq:1.1}, \eqref{eq:1.5} and \eqref{eq:5.1}, for $i=1,\dots,N$, we have 
\begin{equation}
\label{eq:5.3}
\begin{split}
 & \frac{u_{x_i}(x,t)^2}{u(x,t)^2}=\left[-\frac{x_i}{2t}+\frac{1}{2t}X_i\right]^2
=\frac{x_i^2}{4t^2}-\frac{x_i}{2t^2}X_i+\frac{1}{4t^2}X_i^2,\\
 & \frac{u_{x_ix_i}(x,t)}{u(x,t)} =-\frac{1}{2t}+\frac{x_i^2}{4t^2}-\frac{x_i}{2t^2}X_i+\frac{1}{4t^2}Y_i,
\end{split}
\end{equation}
for $(x,t)\in{\bf R}^N\times(0,\infty)$, where 
$$
X_i:=\frac{1}{u}\int_{{\bf R}^N}y_i G(x-y,t)u_0(y)\,dy,\qquad
Y_i:=\frac{1}{u}\int_{{\bf R}^N}y_i^2G(x-y,t)u_0(y)\,dy.
$$
It follows from $\mbox{supp}\,u_0\subset B(0,R)$ that 
\begin{equation}
\label{eq:5.4}
|X_i|\le R,\qquad 0\le Y_i\le R^2.
\end{equation}
Since 
\begin{equation*}
\begin{split}
v_{x_i}= & -\gamma(-\log u)^{-(1-\gamma)}\frac{u_{x_i}}{u},\\
v_{x_ix_i}= & -\gamma(1-\gamma)(-\log u)^{-(2-\gamma)}\frac{(u_{x_i})^2}{u^2}
+\gamma(-\log u)^{-(1-\gamma)}\frac{(u_{x_i})^2}{u^2}
-\gamma(-\log u)^{-(1-\gamma)}\frac{u_{x_ix_i}}{u},
\end{split}
\end{equation*}
by \eqref{eq:5.3} we obtain
\begin{equation*}
\begin{split}
\frac{2t}{\gamma}(-\log u)^{1-\gamma}v_{x_ix_i}
 & =2t\left[-\frac{u_{x_ix_i}}{u}+\frac{(u_{x_i})^2}{u^2}\right]-2t(1-\gamma)(-\log u)^{-1}\frac{(u_{x_i})^2}{u^2}\\
 & =1+\frac{1}{2t}X_i^2-\frac{1}{2t}Y_i
 -(1-\gamma)(-\log u)^{-1}\left[\frac{x_i^2}{2t}-\frac{x_i}{t}X_i+\frac{1}{2t}X_i^2\right]\\
\end{split}
\end{equation*}
for $(x,t)\in{\bf R}^N\times(T,\infty)$.
Since $\lim_{t\to\infty}\|u(t)\|_{L^\infty({\bf R}^N)}=0$ (see \eqref{eq:5.2}), 
taking a sufficiently large $T$ if necessary, 
we have 
\begin{equation}
\label{eq:5.5}
\frac{2t}{\gamma}(-\log u)^{1-\gamma}v_{x_ix_i}
\ge 1+\frac{1}{4t}X_i^2-\frac{1}{2t}Y_i
-(1-\gamma)(-\log u)^{-1}\left[\frac{x_i^2}{2t}-\frac{x_i}{t}X_i\right]
\end{equation} 
for $(x,t)\in{\bf R}^N\times(T,\infty)$.

Let $0<\epsilon<1$. 
By \eqref{eq:5.2} we take a sufficiently large $T$ so that  
\begin{equation}
\label{eq:5.6}
(-\log u)^{-1}\le\left(-\log\biggr[(4\pi t)^{-\frac{N}{2}}\|u_0\|_{L^1({\bf R}^N)}\biggr]\right)^{-1}
\le\left(\frac{N}{4}\log t\right)^{-1}
\end{equation}
for $(x,t)\in{\bf R}^N\times(T,\infty)$. 
We consider the case where $(x,t)\in {\bf R}^N\times(T,\infty)$ with $|x|^2\le\epsilon t\log t$. 
By \eqref{eq:5.4} and \eqref{eq:5.6} we have 
\begin{equation}
\label{eq:5.7}
\begin{split}
(-\log u)^{-1}\left[\frac{x_i^2}{2t}-\frac{x_i}{t}X_i\right]
 & \le\left(\frac{N}{4}\log t\right)^{-1}
\left[\frac{|x|^2}{2t}+R\frac{|x|}{t}\right]\\
 & 
\le\left(\frac{N}{4}\log t\right)^{-1}
\left[\frac{\epsilon \log t}{2}+R\epsilon^{\frac{1}{2}}t^{-\frac{1}{2}}(\log t)^{\frac{1}{2}}\right]\\
 & 
=\frac{2}{N}\epsilon+\frac{4}{N}R\epsilon^{\frac{1}{2}}t^{-\frac{1}{2}}(\log t)^{-\frac{1}{2}}.
\end{split}
\end{equation} 
By \eqref{eq:5.4}, \eqref{eq:5.5} and \eqref{eq:5.7}, 
taking a sufficiently small $\epsilon>0$ and a sufficiently large $T$ if necessary, 
we obtain 
\begin{equation}
\label{eq:5.8}
\frac{2t}{\gamma}(-\log u)^{1-\gamma}v_{x_ix_i}
\ge 1-\frac{R^2}{2t}
 -(1-\gamma)\left[\frac{2}{N}\epsilon+\frac{4}{N}R\epsilon^{\frac{1}{2}}t^{-\frac{1}{2}}(\log t)^{-\frac{1}{2}}\right]
\ge\frac{1}{2}
\end{equation} 
for $(x,t)\in{\bf R}^N\times(T,\infty)$ with $|x|^2\le\epsilon t\log t$. 

We consider the case where $(x,t)\in {\bf R}^N\times(T,\infty)$ with $|x|^2>\epsilon t\log t$. 
Let $\delta$ be a positive constant to be chosen later. 
Since $\mbox{supp}\,u_0\subset B(0,R)$, 
by \eqref{eq:1.5}, taking a sufficiently large $T$ if necessary, 
we have 
$$
u(x,t)\le (4\pi t)^{-\frac{N}{2}}\exp\left(-\frac{|x|^2}{4(1+\delta)t}\right)\|u_0\|_{L^1({\bf R}^N)}.
$$
This implies that 
\begin{equation*}
 (-\log u)^{-1}
 \le \left[\frac{N}{2}\log(4\pi t)+\frac{|x|^2}{4(1+\delta)t}-\log\|u_0\|_{L^1({\bf R}^N)}\right]^{-1}
 \le\frac{4(1+\delta)t}{|x|^2}.
\end{equation*}
It follows from \eqref{eq:5.4} and \eqref{eq:5.5} that 
\begin{equation}
\label{eq:5.9}
\begin{split}
\frac{2t}{\gamma}(-\log u)^{1-\gamma}v_{x_ix_i}
 & \ge 1-\frac{R^2}{2t}
 -(1-\gamma)\frac{4(1+\delta)t}{|x|^2}\left[\frac{|x|^2}{2t}+R\frac{|x|}{t}\right]\\
 & \ge 1-2(1-\gamma)(1+\delta)-\frac{R^2}{2t}-(1-\gamma)R\frac{4(1+\delta)}{(\epsilon t\log t)^{\frac{1}{2}}}.
\end{split}
\end{equation} 
Since $1/2<\gamma\le 1$, taking a sufficiently small $\delta>0$, 
we see that $1-2(1-\gamma)(1+\delta)\ge\delta$. 
Then, by \eqref{eq:5.9}, taking a sufficiently large $T$ if necessary, 
we obtain 
\begin{equation}
\label{eq:5.10}
\frac{2t}{\gamma}(-\log u)^{1-\gamma}v_{x_ix_i}\ge\frac{\delta}{2} 
\end{equation} 
for  $(x,t)\in {\bf R}^N\times(T,\infty)$ with $|x|^2>\epsilon t\log t$. 
Combining \eqref{eq:5.8} and \eqref{eq:5.10},  
we deduce that $v(\cdot,t)$ is convex in ${\bf R}^N$ for $t\ge T$. 
Therefore we see that 
$L_\alpha(u(\cdot,t))$ is concave in ${\bf R}^N$ for $t\ge T$. 
This together with Lemma~\ref{Lemma:2.1}  implies that 
$u(\cdot,t)$ is $\alpha$-logconcave in ${\bf R}^N$ for $t\ge T$. 
Thus Theorem~\ref{Theorem:3.4} follows.
$\Box$
\medskip

\noindent
{\bf Acknowledgements.} 
The first author was supported in part by the Grant-in-Aid for Scientific Research (A)(No.~15H02058)
from Japan Society for the Promotion of Science.
The second author has been partially supported by INdAM through a GNAMPA Project.
The third author was supported in part by the Grant-in-Aid for Young Scientists (B)(No.~15K17536)
and the Grant-in-Aid for Scientific Research (C)(No.~16KT0132).

\end{document}